\documentclass[10pt]{article}
\usepackage{amsmath, amsfonts, amssymb, amsthm, newlfont, graphicx}
\usepackage{epstopdf}

\usepackage{color, pdfcolmk}
 \usepackage[normalem]{ulem}

\newcommand{\be}{\begin{equation}}
\newcommand{\ee}{\end{equation}}

\begin{document}

\title{The generalized Taylor series approach is not equivalent to the homotopy analysis method}

\author{Robert A. Van Gorder$^*$ \\  \small  Mathematical Institute, University of Oxford\\
\small Andrew Wiles Building, Radcliffe Observatory Quarter, Woodstock Road, Oxford, OX2 6GG, United Kingdom\\
\small Email: Robert.VanGorder@maths.ox.ac.uk}       
\date{\today}       
\maketitle

\begin{abstract}
In recent work on the area of approximation methods for the solution of nonlinear differential equations, it has been suggested that the so-called generalized Taylor series approach is equivalent to the homotopy analysis method. In the present paper, we demonstrate that such a view is only valid in very special cases, and in general the homotopy analysis method is far more robust. In particular, the equivalence is only valid when the solution is represented as a power series in the independent variable. As has been shown many times, alternative basis functions can greatly improve the error properties of homotopy solutions, and when the base functions are not polynomials or power functions, we no longer have that the generalized Taylor series approach is equivalent to the homotopy analysis method. We demonstrate this by consideration of an example where the generalizes Taylor series must always have a finite radius of convergence (and hence limited applicability), while the homotopy solution is valid over the entire infinite domain. We then give a second example for which the exact solution is not analytic, and hence it will not agree with the generalized Taylor series over the domain. We conclude that the generalized Taylor series approach is not equivalent to the homotopy analysis method. In particular, the generalized Taylor series can at best recover local information about solutions, whereas the homotopy analysis method can recover the solutions globally if appropriate base functions are selected. Such results have important implications for how series are calculated when approximating solutions to nonlinear differential equations.
\end{abstract}

\section{Introduction}
The Homotopy Analysis Method (HAM) is an analytical solution method which allows one to approximate the solution to nonlinear ordinary differential equations, partial differential equations, integral equations, and so on \cite{ham1,ham0,ham2}. The HAM has proven useful for a variety of such problems \cite{ham3,ham4,ham6,ham6a,ham6b}, owing to the fact that it is unique among analytical or perturbation methods in that it gives a way to minimize the error of approximations by way of an auxiliary parameter, commonly referred to as a convergence control parameter. For instance, one may minimize the error or residual error of approximate solutions over all possible choices of this parameter, and this process is referred to as the optimal HAM (or, OHAM); see \cite{ham7,ham8,ham9,ham10}. The HAM also gives one great freedom in selecting the form of the solutions via representation of solutions \cite{ham1,ham11}, since one has control over the type of basis functions employed in such a representation. 

In a series of papers, Liu \cite{gts1,gts2,gts3} claimed that the HAM was actually equivalent to obtaining a generalized Taylor series expansion at some point in the problem domain. Liu was able to interpret the convergence control parameter in the HAM as relating to the expansion location a generalized Taylor series. A couple of specific examples were provided to demonstrate these points, all for ODEs with rather well-behaved solutions. In addition to those papers, the HAM has been applied in a number of cases in ways that make it seem equivalent to a Taylor series approach (we do not list all such instances here).

The aim of this paper is to show that although the observations of Liu do hold in some specific examples, it is not in general true that the generalized Taylor series is equivalent to the HAM solution. Rather, this equivalence can only hold when the HAM solution is represented in a power series of a single independent variable. In other situations, such as when the HAM series is represented in some other non-polynomial function basis, the HAM can give series solutions which converge everywhere, even though no generalized Taylor series can converge over the whole problem domain. These points and clarifications are strongly worth making, in light of the fact that the work of Liu has been cited in other works, and hence we would prefer to cast light on the differences in the generalized Taylor series and HAM approaches.

In Section 2 we discuss the connection between HAM and the generalized Taylor series method. In Section 3, we consider solutions to a nonlinear boundary value problem which admits at best a local representation via any generalized Taylor series, yet a global representation via HAM. In Section 4, we consider solutions to a nonlinear boundary value problem which is not analytic at the origin. For this problem, the solution is not equal to it's power series at the origin. However, using a different set of basis functions (which themselves are not analytic at the origin) we are able to construct the global exact solution using HAM.

\section{Mathematical underpinnings of the generalized Taylor series}
In HAM, one attempts to solve a nonlinear (ordinary or partial) differential equation $N[u]=0$ by solving a related problem that is more amenable to analysis. Introducing an auxiliary linear operator, $L$, and convergence control parameter, $h$, and one then constructs the homotopy of operators 
\be 
H[u] = (1-q)L[u] - hqN[u]\,,
\ee
which is itself a differential operator. Here $q\in [0,1]$ is an embedding parameter, such that upon setting $q=0$ we have $L[u]=0$ implies $H[u]=0$, while setting $q=1$ gives $N[u]=0$ implies $H[u]=0$. Setting $H[u]\equiv 0$ identically, we have a homotopy of the operators $L$ and $N$. If we then consider the more general equation $H[u]=0$, then if a solution exists it should depend on the embedding parameter $q$, say $u(t) = \hat{u}(t;q)$. We should have that $\hat{u}(t;0)$ is a solution of $L[u] =0$, while $\hat{u}(t;1)$ is a solution is $N[u]=0$. If $\hat{u}(t;q)$ varies continuously in $q\in [0,1]$, then taking $q$ from 0 to 1 gives a map from a solution of $L[u]=0$ (a linear problem) to a solution of $N[u]$ (the original nonlinear problem of interest). Note that an entire literature exists for the case where $h = -1$, in which case the method is often referred to as the homotopy perturbation method (HPM); see \cite{hpm1,hpm2,hpm3}. On the other hand, when $h$ is picked in a specific way to minimize the error in the approximate solution formed through truncation, the method is often referred to as the optimal homotopy analysis method (OHAM); see \cite{ham7,ham8,ham9,ham10}.

We shall not go into details of the HAM solution procedure for the sake of brevity; see \cite{ham1,ham2,ham3,ham11} for details and many worked examples. Still, we do need to remark that one has great freedom to pick the operator $L$ in HAM, as has been previously discussed \cite{ham1,ham2,ham3,ham11}. Importantly, it has been shown that the selection of $L$ can result in a specific set of base functions used in the solution representation. In contrast, in perturbation methods or other approaches where the linearized equations are of a fixed form, the base function in the solution representation are fixed. If $L$ is picked so that a power series representation for a solution is obtained, then we effectively obtain a power series solution. When a function agrees with its power series, then that series representation is unique, and therefore the series representation is equal to the Taylor series for that function. 

As a corollary to this, we know that if a power series solution in a variable is found via HAM, then that series must be equivalent to the generalized Taylor series with all terms calculated by centering the series at some fixed point in the domain \cite{gts1}. In particular, for some initial value problem involving an unknown function $u(t)$, if one picks a HAM series solution in powers of the independent variable $t$, one obtains the series solution
\be
u(t) = \sum_{n=0}^\infty A_n(h)t^n\,.
\ee
The framework of the generalized Taylor series method outlined in \cite{gts1} tells us that there exists a corresponding and equivalent series representation (once all terms are multiplied through) of the form
\be\label{gt}
u(t) = \sum_{n=0}^\infty a_n (t-t_0(h))^n  = \sum_{n=0}^\infty \frac{f^{(n)}(t_0(h))}{n!} (t-t_0(h))^n\,,
\ee
where the resulting series is expanded about some $h$-dependent point of expansion, $t_0(h)$, rather than zero. This is a more concise way to represent the power series solution, since $h$ enters in exactly one way, as the point the series is expanded about. See \cite{gts2,gts3} for proofs. Therefore, when the HAM solution takes the form of a power series, we may express the solution more concisely via the generalized Taylor theorem.

Note that the equivalence between is only true when the HAM gives a power series solution. If a solution is given in terms of a series of some other functions, the result no longer holds. Indeed, if basis functions can be chosen that have better global properties than power series, then we may be able to enhance the region of convergence of the solutions. A relevant example of this is given in the following section.

To apply the HAM, we assume solutions which are analytic in the embedding parameter, $q$. However, we do not require that solutions need to be analytic in the independent variable(s) of the problem. Yet, if one is to attempt a generalized Taylor series solution, it is clear that such a series can only agree with the true solution of a problem provided that that solution is analytic on the problem domain. Intuitively, then, the generalized Taylor series formulation can recover a HAM solution only in some special cases, and in general cannot be used to find that HAM solution. We provide two examples to illustrate our point.

\section{An example for which the generalized Taylor series is valid only locally}
Consider the boundary value problem
\be \label{ode}
u'' + 2u^3 -u =0\,,\quad  u(0) = 1\,, \quad \text{and} \quad \lim_{|t|\rightarrow \infty}u(t) = 0\,.
\ee
The exact solution of the boundary value problem \eqref{ode} is given by $u(t) = \text{sech}(t)$.

In light of the boundary conditions \eqref{ode}, it makes sense to consider base functions that decay as $t\rightarrow \infty$. We shall pick an auxiliary linear operator $L$ such that we obtain base functions $e^{-t}, e^{-2t}, e^{-3t},\dots$. Such a choice of $L$ is given by $L [u ] = u'' - u$. The inversion $L^{-1}[0] = c_1 e^{t} + c_2 e^{-t}$ gives the general form of the order zero solution. Using the boundary conditions \eqref{ode}, we have $c_1 =0$ and $c_2 =1$. Hence, $U_0(t) = e^{-t}$. Calculating the higher order terms in the HAM solution expansion, one may shown that
\be 
U(t;h) = U_0(t) + \sum_{n=1}^\infty U_n(t;h)  = \sum_{n=1}^\infty \nu_n(h)e^{-nt}  \,,
\ee
where each $\nu_n(h)$ is some function of $h$, for $n=1,2,3,\dots$. One may calculate the $\nu_n(h)$ to obtain
\be 
\nu_n(-1) = \begin{cases}
2(-1)^m & \text{if} \quad n=2m+1\,,\\
0 & \text{if}\quad n = 2m\,.
\end{cases}
\ee
This gives the representation
\be 
U(t;-1) = \sum_{m=0}^\infty 2(-1)^m e^{-(2m+1)t}\,.
\ee
This series is convergent over the whole interior of the domain, $t\in (0,\infty)$.

Note that the true solution of this boundary value problem is $u(t) = \text{sech}(t)$. However, it is well known that $\text{sech}(t)$ has a finite radius of convergence, and the power series representation for $\text{sech}(t)$, which is given by
\be 
\text{sech}(t) = 1 - \frac{1}{2}t^2 + \frac{5}{24}t^4 - \frac{61}{720}t^6 + \cdots + \frac{(-1)^k E_k}{(2k)!}t^{2k} + \cdots
\ee
is valid for $|t|<\frac{\pi}{2}$. This series is therefore valid if we take the boundary condition at $t=0$, but we cannot take into account the condition as $t\rightarrow \infty$. If we center the series elsewhere, we encounter similar issues (we would shift the region of convergence, yet it would still consist of a finite bounded interval). To see why, note that while many nonlinear ordinary differential equations arising in physical applications result in real-valued solutions, it is well known from the analytic theory of such equations that a series solution's region of convergence is restricted by the appearance of poles in the complex plane, not just on the real line. This is important to note, as many solutions to nonlinear differential equations can have complex poles, even if the behavior on the real line is sufficiently nice. In our example, note that while $\text{sech}(t)$ is bounded on the real line, it does have poles in the complex plane. The nearest poles to the real axis are $t = \pm \frac{i\pi}{2}$. It is these poles that give the bound on the region of convergence of $|t|<\frac{\pi}{2}$ when we construct a Taylor series solution at the origin. Meanwhile, if we attempt to construct a Taylor series solution centered at some other point $t=t_0 >0$, then the distance to the nearest pole in the complex plane is $d = \sqrt{t_0^2 + \frac{\pi^2}{4}}$. Therefore, a Taylor series solution for $\text{sech}(t)$ centered at $t=t_0>0$ will converge on the finite region
\be 
t_0 - \sqrt{t_0^2 + \frac{\pi^2}{4}} < t < t_0 + \sqrt{t_0^2 + \frac{\pi^2}{4}}\,.
\ee
When $t_0 =0$, this gives $|t|<\frac{\pi}{2}$, as expected. Therefore, no matter the choice of $t_0 >0$, the convergence region for a Taylor series solution to $\text{sech}(t)$ is always finite. Hence, a Taylor series cannot be the solution over the boundary value problem \eqref{ode} over the whole domain $t\in (0,\infty)$. This highlights the need to search for alternate solution forms if we are to obtain convergent solution over the whole problem domain. Therefore, a Taylor series (or, a generalized Taylor series) cannot give a solution to the boundary value problem \eqref{ode} over the domain $t\in (0,\infty)$. Rather, it can only give a local approximation.

It was shown that the generalized Taylor series solution recovers exactly the HAM solution centered at $t=0$. We should remark that this HAM solution was calculated under the assumption of a power series solution expression. When applying the HAM, one has great freedom in selecting the auxiliary linear operator, which in turn gives a specific set of basis function. Therefore, by selecting alternative auxiliary linear operators, we can obtain HAM solutions for this problem which are actually valid over the entire problem domain. To do this, we seek an auxiliary linear operator which gives basis functions that are uniformly bounded over the semi-infinite problem domain $[0,\infty)$. Of course, powers of the independent variable will become unbounded as we tend toward the far-field, making it impossible to accurately enforce the far-field condition. This is why we selected the linear operator $L$ as we did.

Note that we can directly obtain the series solution in the basis of functions $e^{-nt}$ directly, if we know that the solution of the boundary value problem \eqref{ode} is $u(t)=\text{sech}(t)$. Indeed, by definition,
\be 
\text{sech}(t) = \frac{2e^{-t}}{1+e^{-2t}} = 2e^{-t}\sum_{m=0}^\infty (-e^{-2t})^m = \sum_{m=0}^\infty 2(-1)^m e^{-(2m+1)t}\,,
\ee
where the manipulations are valid on the interior of the domain, $t\in (0,\infty)$. Note that for $t<0$ one can also obtain a series representation,
\be 
\text{sech}(t) = \frac{2e^{t}}{1+e^{2t}} = 2e^{t}\sum_{m=0}^\infty (-e^{2t})^m = \sum_{m=0}^\infty 2(-1)^m e^{(2m+1)t}\,,
\ee
valid on $t\in (-\infty,0)$. Therefore, we have the global representation
\be 
\text{sech}(t) =\begin{cases}
\sum_{m=0}^\infty 2(-1)^m e^{(2m+1)t} & \text{if} \quad t<0\,,\\
1 & \text{if} \quad t=0\,,\\
\sum_{m=0}^\infty 2(-1)^m e^{-(2m+1)t}& \text{if} \quad t>0\,,
\end{cases}
\ee
and this representation is valid for all $t\in \mathbb{R}$. This is exactly the HAM solution for $h=-1$.

\begin{figure}
\begin{center}
\includegraphics[scale=0.4]{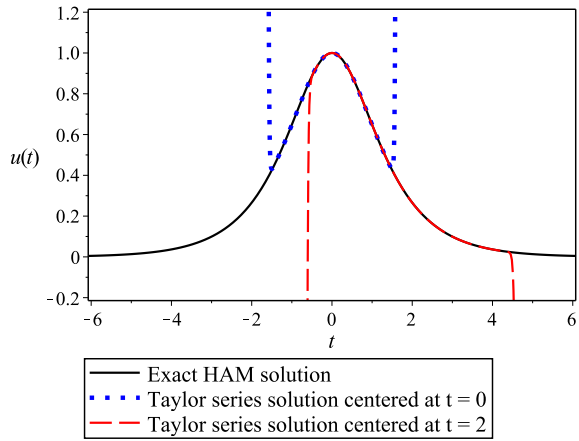}
\end{center}
\vspace{-0.3in}
\caption{Plot of the exact solution to \eqref{ode} (which coincides with the HAM solution when the auxiliary operator $L$ is employed) and two Taylor series solutions, centered at $t=0$ and $t=2$, respectively. The order of these series is $O(t^{100})$, hence we are very confident of their accuracy over the regions on which they converge. As we see, adjusting the expansion point for the generalized Taylor series can modify the region of convergence, although this region will always be finite. In contrast, using non-power series, the HAM solution is able to reproduce the global solution on $\mathbb{R}$.\label{Fig1}}
\end{figure}

In Figure 1, we plot the exact solution obtained via HAM with $h=-1$ and $L$ as given above. We also plot the high order Taylor series solutions centered at $t=0$ and $t=2$ (corresponding to different values of $h$ if one assumes a HAM solution in terms of polynomial base functions). To get a good understanding of the behavior of the Taylor series, we include up to order $t^{100}$ terms in our calculations. We see that while changing the location $t=t_0$ at which the Taylor series is centered does expand the region of convergence, this region is still finite for any finite $t_0$. Therefore, no generalized Taylor series can be constructed which converges over the whole domain $t\in (0,\infty)$.

\section{An example with non-analytic solutions in the independent variable}
Consider the boundary value problem 
\be \label{ode2}
uu'' - {u'}^2 + 2t^{-3}u^2 =0\,,\quad u(-\infty) = 0\,, \quad \text{and} \quad u(\infty) = 1\,.
\ee
From the form of \eqref{ode2} we see that this equation is non-analytic. Yet, this equation has the exact solution
\be \label{sol2}
u(t) = \begin{cases}
\exp\left( -t^{-1} \right) & \text{for} ~~ t>0\,,\\
0 & \text{for} ~~ t\leq 0\,.
\end{cases}
\ee
Note that this solution is not only continuous, but also smooth. However, one may show that this solution will not agree with its Taylor series at $t=0$, hence it is not analytic. (Indeed, $\exp(-t^{-t})$ is a standard example of such a function.) As such, no Taylor series in $t$ of the form \eqref{gt} can represent this solution of \eqref{ode2}, so the generalized Taylor series approach cannot be applied to this problem.

\begin{figure}
\begin{center}
\includegraphics[scale=0.4]{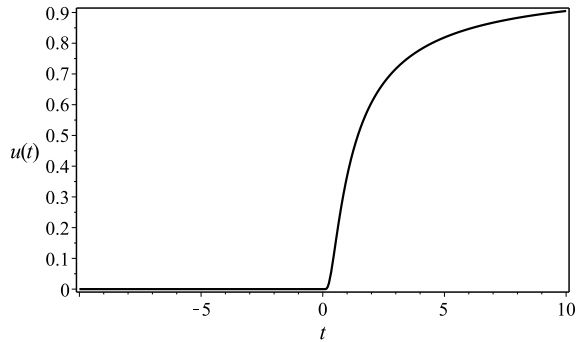}
\end{center}
\vspace{-0.3in}
\caption{Plot of the exact solution \eqref{sol2} to \eqref{ode2} (which coincides with the HAM solution when the auxiliary operator $L_2$ is employed). Although the solution is not analytic at $t=0$, it can be reproduced by HAM using non-analytic base functions. \label{Fig2}}
\end{figure}

With that said, one can search for a solution via HAM in some other base functions which are not simply powers of $t$. Consider the linear operator $L_2[u]$ given by $L_2[u] = u'' - t^{-2}u' + 2t^{-3}u$. A homogeneous solution of $L_2[u] = 0$ is given by
\be 
u(t) = C_1 \exp\left( - t^{-1}\right) + C_2\left\lbrace t + \exp\left( - t^{-1}\right)\text{Ei}\left(1, -t^{-1}\right) \right\rbrace\,,
\ee
where Ei denotes the relevant exponential integral. Picking $C_1=1$ and $C_2=0$ gives the solution for $t>0$, while picking $C_1 = C_2 =0$ gives the solution for $t\leq 0$. Matching these two solution branches, we obtain the solution valid over the real line. Performing HAM with this choice of auxiliary linear operator will then allow us to recover the exact solution \eqref{sol2}. For this example, we have constructed the linear operator to give the solution on the zeroth order approximation, but similar comments would follow if one needed to calculate higher order approximations to construct a solution. We give a plot of this solution in Figure 2.

Note that we are able to recover the solution \eqref{sol2} by choosing an auxiliary linear operator which is not analytic. This is fine, as technically the homotopy need only be analytic in the embedding parameter, $q$. In this way, one may use the HAM to obtain a solution which satisfies an ODE which is not analytic in the independent variable(s) (yet which still has an exact solution). 

\section{Discussion}

In the optimal HAM, one attempts to minimize some measure of the error of the approximate analytical solution by use of the convergence control parameter, $h$. Although one can use the generalized Taylor series approach to modify the interval of convergence, the approach is not very useful if relatively few terms are taken and then one tries to minimize error. For most interesting applications, only relatively few terms can be calculated due to the strong nonlinearity inherent in many real world problems. For this reason, the generalized Taylor series alone does not offer many benefits over the traditional Taylor series solutions for differential equations. It is well-known that for strongly nonlinear problem, often have finite domains of convergence and hence limited regions of applicability. For such cases, non-polynomial base functions will be far more useful.

What this suggests is that the role of the convergence control parameter, $h$, in the HAM is not as clear cut as is suggested by \cite{gts1}. Furthermore, we have shown that the auxiliary linear operator can be selected to overcome shortcomings one would encounter by simply selecting a Taylor series solution. Indeed, this choice is strongly related to the rule of solution expression, which permits one to construct solutions in non-polynomial base functions, resulting in solutions which are not simple power series. We conclude that the generalized Taylor series approach is not equivalent to the homotopy analysis method.

We should also remark that power series are far less useful for partial differential equations, since it is far harder to control the convergence of multivariate power series \cite{power}, and indeed any applications of the so-called generalized Taylor series have been limited to ordinary differential equations. It is also worth noting that power series for multivariate functions may have complicated convergence regions. On the other hand, the HAM allows one to pick other basis functions that can be required to satisfy sufficient boundedness conditions over a multiple dimension domain \cite{ham1,ham11}. This illustrates one large advantage that the HAM has over the generalized Taylor series approach when the domain of the problem is of dimension greater than one.

\clearpage

\end{document}